\documentclass[journal]{IEEEtran} %option lettersize useless?
\usepackage{amsmath,amsfonts}
\usepackage{algorithmic}
\usepackage{algorithm}
\usepackage{array}
\usepackage[caption=false,font=normalsize,labelfont=sf,textfont=sf]{subfig}
\usepackage{textcomp}
\usepackage{stfloats}
\usepackage{url}
\usepackage{verbatim}
\usepackage{graphicx}
\usepackage{cite}
\usepackage{hyperref}
\usepackage{comment}
\usepackage{siunitx}
\usepackage{soul}

\newcommand\bs[1]{\boldsymbol{#1}} %vector
\newcommand\mbf[1]{\mathbf{#1}} %algebraic vector
\newcommand\mbb[1]{\mathbb{#1}} %auto math \mathbb 
\newcommand\mcal[1]{\mathcal{#1}}

\newcommand\J{\bs{J}}
\newcommand\M{\bs{M}}
\newcommand\Jd{\mbf{j}}
\newcommand\Md{\mbf{m}}

\newcommand\nh{\bs{\hat{n}_r}}
\newcommand\Tc{\mcal{T}}
\newcommand\Kc{\mcal{K}}
\newcommand\Ic{\mcal{I}}
\newcommand\Pmc{\mcal{P}^-}
\newcommand\Rc{\mcal{R}}
\newcommand\Trwg{\mbf{T}}
\newcommand\Krwg{\mbf{K}}

\newcommand\Tbc{\mbb{T}}
\newcommand\Kbc{\mbb{K}}
\newcommand\Gbc{\mbb{G}}
\newcommand\Mrwg{\mbf{P}_k}
\newcommand\Mbc{\mbb{P}_k}
\newcommand\Pl{\mbb{P}^\Lambda}
\newcommand\Ps{\mbf{P}^\Sigma}
\newcommand\Psh{\mbb{P}^{\Sigma H}}
\newcommand\Plh{\mbf{P}^{\Lambda H}}
\newcommand\Lam{\mbf{\Lambda}}
\newcommand\Sig{\mbf{\Sigma}}

\newcommand\lrr[1]{\left(#1\right)}
\newcommand\lrs[1]{\left[#1\right]}

\newcommand\lra[1]{\langle#1\rangle}
\newcommand\diff{\,\mathrm{d}}

\newcommand\im{\text{i}}

\newcommand\cald{Calder\'{o}n}
\newcommand\poin{Poincaré}

\usepackage{lipsum} 
\usepackage[dvipsnames]{xcolor}
\usepackage{todonotes}
\begin{document}

%\title{On the zero-field constraint of a Poincaré-Steklov-based single source method}
\title{\begin{color}{black}Stabilized Single Current Inverse Source Formulations Based on Steklov-\poin{} Mappings\end{color}}

\author{IEEE Publication Technology Department
\thanks{Manuscript created October, 2020; This work was developed by the IEEE Publication Technology Department. This work is distributed under the \LaTeX \ Project Public License (LPPL) ( http://www.latex-project.org/ ) version 1.3. A copy of the LPPL, version 1.3, is included in the base \LaTeX \ documentation of all distributions of \LaTeX \ released 2003/12/01 or later. The opinions expressed here are entirely that of the author. No warranty is expressed or implied. User assumes all risk.}}

\author{Paolo Ricci, Ermanno Citraro, Adrien Merlini, Francesco P. Andriulli}

\author{Paolo~Ricci~\IEEEmembership{Graduate Student Member, IEEE}, Ermanno~Citraro~\IEEEmembership{Graduate Student Member, IEEE}, Adrien~Merlini~\IEEEmembership{Member, IEEE}, and Francesco~P.~Andriulli~\IEEEmembership{Senior Member, IEEE}
\thanks{Paolo~Ricci is with the Department of Electronics and Telecommunications, Politecnico di Torino, Turin, Italy.}
\thanks{Ermanno~Citraro is with the Department of Electronics and Telecommunications, Politecnico di Torino, Turin, Italy.}
\thanks{Adrien~Merlini is with the Microwave Department, IMT Atlantique, Brest, France.}
\thanks{Francesco~P.~Andriulli is with the Department of Electronics and Telecommunications, Politecnico di Torino, Turin, Italy.}}

\markboth{ Journal of \LaTeX\ Class Files,~Vol.~18, No.~9, September~2020}%
{How to Use the IEEEtran \LaTeX \ Templates}


\maketitle

\begin{abstract}
 The inverse source problem in electromagnetics has proved quite relevant for a large class of applications. In antenna diagnostics in particular, Love solutions are often sought at the cost of an increase of the dimension of the linear system to be solved. In this work, instead, we present a reduced-in-size single current formulation of the inverse source problem that obtains one of the Love currents via a stable discretization of the Steklov-\poin{} boundary operator leveraging dual functions. The new approach is enriched by theoretical treatments and by a further low-frequency stabilization \textcolor{black}{of the Steklov-\poin{} operator} based on the quasi-Helmholtz projectors that is the first of its kind in this field. The effectiveness and practical relevance of the new schemes are demonstrated via both theoretical and numerical results.
\end{abstract}

\begin{IEEEkeywords}
Inverse source problem, boundary element method, Steklov-\poin{} operator, Love currents, low-frequency breakdown.
\end{IEEEkeywords}

\section{Introduction}

\IEEEPARstart{T}{he} inverse \textcolor{black}{ source problem in electromagnetics, i.e. the recovery of a configuration of sources radiating a given measured field, has been adopted in a variety of applications ranging from antenna diagnostics to near-to-far-field reconstructions\cite{lopez_improved_2009, cappellin_advanced_2011, foged_practical_2012}.}
These sources are often electric and/or magnetic current distributions residing on a conveniently placed equivalent surface that can be tailored to scatter the target field by virtue of the equivalence theorem. These currents have traditionally been found within a boundary element framework on apertures or on arbitrary equivalent surfaces (see for example \cite{petre_planar_1992,sarkar_near-field_1999}).
%More recently, a general approach which takes into account both types of currents has been considered and its performance studied against the former so-called single current formulations.
%Single current solutions, that reconstruct either the electric or magnetic current, are appealing because of the reduced dimensions of the linear system to be solved, but they been reported to require more care than their dual current counterparts \cite{kornprobst_solution_2019, quijano_field_2010}.
\textcolor{black}{Among inverse source strategies, single current solutions, that reconstruct only one family among electric or magnetic currents, are appealing because of the reduced dimensions of the linear systems to be solved and because of their reduced (numerical) nullspace that is limited to the intrinsic ill-posedness of the problem associated to the non-radiating modes. These strategies, however, have been reported to require more care in the solution process if further physical constraints are not used to ensure a simple relationship between equivalent currents and fields \cite{quijano_field_2010, kornprobst_solution_2019}}.
%Upon discretization, however, both families of approaches will exhibit numerical nullspaces corresponding to t
On the other hand, the operators involved in double current formulations exhibit an additional, large nullspace corresponding to the excess in discrete degrees of freedom; the solution of the corresponding linear systems require pseudoinversions that have been efficiently performed via iterative solvers \cite{alvarez_reconstruction_2007, eibert_multilevel_2009, eibert_inverse_2010}.
While the nullspace of the system can partially be handled in such a way, the numerical ill-conditioning of the matrix, inherited by the ill-posed nature of the inverse problem, remains to be addressed.
To this end, truncated singular value decompositions (TSVD) or Tikhonov regularizations have been used to further regularize the problem \cite{jorgensen_improved_2010, cappellin_advanced_2011}.

Another feature of interest among inverse source schemes is their capacity to find equivalent Love currents---that are directly related to the tangential fields---which is considered in literature particularly useful for antenna diagnostics \cite{quijano_field_2010, jorgensen_improved_2010}.
The Love currents can be obtained by adding further constraints to double current formulations \cite{quijano_field_2010, kornprobst_accuracy_2021, phaneuf_formulation_2021} or by filtering any of the solution via \cald{} projection \cite{kornprobst_inverse_2019}.
Another interesting approach, leveraging Huygens radiators and valid for plane waves, has been proposed in \cite{eibert_fast_2016} to reduce the size of the Love-constrained problem to that of a single current formulation, at the price of an approximation. 

\textcolor{black}{In this work we follow a different approach. While still targeting a single current formulation, we leveraged dual discretizations to avoid approximating the relationships linking electric and magnetic currents. The contribution of this paper is then twofold: we present a new single current formulation capable of obtaining Love currents by leveraging a stable discretization of the Steklov-\poin{} operator \cite{de_la_bourdonnaye_formulations_1995} without resorting to any approximations of the electromagnetic relations. This results in a single current formulation that delivers one of the Love currents. Moreover we present the first frequency stabilization 
of Steklov-\poin{} operators via quasi-Helmholtz projectors and we leverage on this new result to stabilize in frequency the new formulations.
%of 
%an analysis establishing that the use of Steklov-\poin{} operators does not prevent the low-frequency stabilization by using quasi-Helmholtz projectors. 
What we propose is then, to the best of our knowledge, the first low-frequency regularization of a full-wave inverse source scheme showing high level of accuracy and numerical stability till arbitrarily low-frequencies.}

%This new formulation has the advantage of a reduced system size typically found in single source formulations . It can be used to obtain magnetic or electric currents, however, when the former are considered, the resulting operator suffers from a low-frequency breakdown inherited from the electric-field integral operator \cite{vecchi_loop-star_1999, jun-sheng_zhao_integral_2000}. To address this issue, we introduce a preconditioner for the system, based on quasi-Helmholtz projectors \cite{andriulli_well-conditioned_2013}, that renders the scheme stable at low frequencies. An early form of the scheme, deprived of theoretical apparatus, has been presented in the conference contribution \cite{ricci_frequency-stabilized_2022}.

The paper is organized as follows: the main electromagnetic operators are introduced in Section~\ref{sec:backround}, the new formulations are presented in Section~\ref{sec:steklov}, while Section~\ref{sec:helmoltz} presents the frequency stabilization of the Steklov-\poin{} operator and its application to the the new equations. Finally Section~\ref{sec:results} illustrates the accuracy and stability of the new formulation through numerical test cases. Section~\ref{sec:conclusion} presents conclusions. Very preliminary results from this work were presented in the conference contribution \cite{ricci_frequency-stabilized_2022}.

% \IEEEpubidadjcol

\section{Background and Notation}
\label{sec:backround}
% \begin{itemize}
%     \item formal description of the ISP and Love currents;
%     \item Steklov Poincaré
%     \item Calder\'on projectors
%     \item continuous operators
%     \item Radiation matrix
%     \item Definition of the (2/3) main strategies in the literature (Vecchi DEqF, korn eta nx mapping, korn con proiettore). Formula only for the mapping (inline)
%     \item ...
% \end{itemize}

Let $\Gamma$ be a smooth manifold in $\mbb{R}^3$ delimiting the internal and external domains $\Omega^-$ and $\Omega^+$. Consider a time-harmonic source in $\Omega^-$ generating Maxwellian fields in $\Omega^-\cup\overline{\Omega^+}$. In light of the equivalence theorem, there exist equivalent current densities $\M$ and $\J$ on $\Gamma$ which radiate in $\Omega^+$ the same fields as the original source and radiate in $\Omega^-$ possibly different electric and magnetic fields; these currents satisfy
%.One possible definition of these currents is
\begin{align}
    \M& = \lrr{\bs{E^+}-\bs{E^{\prime-}}}\label{eq:MBC}\times\nh\,,\\
    \J& = \nh\times\lrr{\bs{H^+}-\bs{H^{\prime-}}}\label{eq:JBC}\,,
\end{align}
where $\nh$ is the unit normal vector to $\Gamma$ in $\bs{r}$ pointing towards $\Omega^+$, $\bs{E}^+$, $\bs{H}^+$ are the original electric and magnetic field in $\Omega^+$ and $\bs{E}^{\prime-}$, $\bs{H}^{\prime-}$ are the new fields in $\Omega^-$. The time-harmonic dependence is assumed and suppressed throughout the paper.

% and $\bs{E}^-$, $\bs{H}^-$ (resp. $\bs{E}^+$, $\bs{H}^+$) are the electric and magnetic fields in $\Omega^-$ (resp. $\Omega^+$).

Solving the inverse source problem consists in finding a set of equivalent currents $\M$, $\J$ given the electric and/or magnetic fields' measurements on a smooth and simply connected manifold $\Gamma_m\subset\Omega^+$. The problem can be solved naturally by the boundary element method. In this framework, define the electric field integral equation (EFIE) operator on $\Gamma$
%\begin{equation}\label{eq:Tdef}
    $\Tc_{\bs{r}}\bs{f} = \im k\,\Tc_{s,\bs{r}}\bs{f}+\im k^{-1}\,\Tc_{h,\bs{r}}\bs{f}$
%\end{equation}
with
% \begin{align}\label{eq:Tshdef}
    $\Tc_{s,\bs{r}}\bs{f} = \nh\times\int_{\Gamma}\frac{e^{\im
    		k\left|\bs{r}-\bs{r}'\right|}}{4\pi\left|\bs{r}-\bs{r}'\right|}{\bs{f}(\bs{r}')}\diff
    \bs{r}'$, $\Tc_{h,\bs{r}}\bs{f} = \nh\times\nabla\int_{\Gamma}\frac{e^{\im
    		k\left|\bs{r}-\bs{r}'\right|}}{4\pi\left|\bs{r}-\bs{r}'\right|}{\nabla_s\cdot\bs{f}(\bs{r}')}\diff
    \bs{r}'$
%\end{align}
and the magnetic field integral equation (MFIE) operator
%\begin{equation}\label{eq:Kdef}
    $\Kc_{\bs{r}}\bs{f}  =  - \nh \times p.v. \int_\Gamma\nabla\times\frac{e^{\im k\left|\bs{r}-\bs{r}'\right|}}{4\pi\left|\bs{r}-\bs{r}'\right|}\bs{f}(\bs{r}')\diff\bs{r}'$,
%\end{equation}
$k$ being the wavenumber, $\bs{r}\in\overline{\Omega^+}$. In the case $\bs{r}\in\Gamma$ $\Tc_{\bs{r}}$, $\Kc_{\bs{r}}$ are denoted by $\Tc$, $\Kc$ respectively. The radiation operator
\begin{equation}\label{eq:Rdef}
    \Rc=\begin{bmatrix}
        -\Kc_{\bs{r}} & \Tc_{\bs{r}}\\
        -\Tc_{\bs{r}} &-\Kc_{\bs{r}}
    \end{bmatrix},
\end{equation}
is a linear map between equivalent sources on $\Gamma$ and measured tangential fields on $\Gamma_m$, meaning that
\begin{equation}\label{eq:linsys}
    \Rc
    \begin{bmatrix}
    -\M\\\eta\J
    \end{bmatrix}
    =
    \begin{bmatrix}
    \nh\times\bs{E}^+\\
    \eta\nh\times\bs{H}^+
    \end{bmatrix}\,,
\end{equation}
with $\eta=\sqrt{\mu/\epsilon}$ and $\epsilon$, $\mu$ being the permittivity and the permeability of the medium respectively.
The inverse problem aims at finding unknown current distributions that satisfy \eqref{eq:linsys}, or part of it. Indeed, by selecting a single block of $\Rc$---either $\Kc_{\bs{r}}$ or $\Tc_{\bs{r}}$---and solving for the corresponding reduced right hand side---$\bs{E}^+$ or $\bs{H}^+$---four different single current formulations can be obtained. Alternatively, three double current formulations can be derived by considering the full radiator or one of its rows only. The latter systems of continuous equations admit several solutions because multiple equivalent currents can radiate the same external field in $\Omega^{+}$ and the physical meaning of the solution depends on the type of implicit or explicit additional constraints used to select a particular solution. The Love currents $\bs{M}_L$, $\bs{J}_L$ are one of these particular solutions that are obtained by imposing the fields radiated in $\Omega^-$ to be identically $\bs{0}$.
%enforcing the constraint that $\bs{E}^-=\bs{0}$, $\bs{H}^-=\bs{0}$ in \eqref{eq:MBC}, \eqref{eq:JBC} respectively.
%After defining the internal \cald{} projector as
One way of enforcing this condition is to leverage the well-known \cald{} projector
\begin{equation}
    \Pmc =
    \begin{bmatrix}
        \frac{\Ic}{2}+\Kc & -\Tc\\
        \Tc & \frac{\Ic}{2}+\Kc
    \end{bmatrix},
\end{equation}
where $\Ic$ is the identity operator, that can be added to the system of equations \eqref{eq:linsys} \cite{kornprobst_accuracy_2021} as
\begin{equation}
    \begin{bmatrix}
    \Rc\\
    \Pmc
    \end{bmatrix}
    \cdot
    \begin{bmatrix}
        \bs{-M}_L\\
        \eta\bs{J}_L
    \end{bmatrix}
    =
    \begin{bmatrix}
        \nh\times\bs{E^+},\,
        \nh\times\eta\bs{H^+},\,
        \bs{0},\,
        \bs{0}
    \end{bmatrix}^T\,.
\end{equation}
%where only one row of $\Rc$ and one of $\Pmc$ are needed to get a non singular matrix.
%In fact, when considering Love currents,
This allows the nullspace of $\Rc$ ($\mcal{N}(\Rc)$), i.e. any solution of the internal problem, to be mapped into itself by $\Pmc$:
\begin{equation}
    \Pmc
    \begin{bmatrix}
        \bs{-M}_L\\
        \eta\bs{J}_L
    \end{bmatrix}
    =
    \begin{cases}
        \begin{bmatrix}\bs{-M}_L\\\eta\bs{J}_L
        \end{bmatrix} & \text{if} 
        \begin{bmatrix}\bs{-M}_L\\\eta\bs{J}_L
        \end{bmatrix}\in\mcal{N}\left(\Rc\right)
        \\
        \bs{0} & \text{otherwise}
    \end{cases}.
\end{equation}

\section{Conforming Discretization of a Steklov-\poin{}-based Equation}
\label{sec:steklov}

In this section we introduce a single source method which enforces the Love condition without increasing the matrix system size with regards to standard single source formulations.
First, consider the Love condition expressed with the inner \cald{} projector
\begin{equation}\label{eq:Pm_sys}
    \Pmc
    \begin{bmatrix}
        -\bs{M}_L\\
        \eta\bs{J}_L
    \end{bmatrix}
    =\bs{0}\,.
\end{equation}
Clearly, for $k$ different from resonant wavenumbers of $\Gamma$, \eqref{eq:Pm_sys} defines a relation between the two Love currents
\begin{equation}
    \eta\bs{J_L} = - \lrr{\frac{\Ic}{2} + \Kc}^{-1} \Tc\lrr{-\bs{M}_L}
    \label{eq:steklov}
\end{equation}
where $\lrr{\frac{\Ic}{2} + \Kc}^{-1} \Tc$ is the Steklov-\poin{} operator \cite{de_la_bourdonnaye_formulations_1995}.
By replacing \eqref{eq:steklov} in the first row equation of \eqref{eq:linsys}, we obtain the \begin{color}{black}equation\end{color}
\begin{equation}
    \lrr{-\Kc_{\bs{r}} - \Tc_{\bs{r}}\lrr{\frac{\Ic}{2} + \Kc}^{-1} \Tc} \lrr{\bs{-M}_L} = \nh\times\bs{E^+}\,,
    \label{eq:radStP}
\end{equation}
which is a single source equation that naturally yields the magnetic Love current $\bs{M}_L$.
  \begin{color}{black}
 %Here we have selected the second row (unknown  $\bs{M}_L$) but 
 If instead of this current, the electric Love current $\bs{J}_L$ is desired as the first outcome of the procedure, a similar strategy can be applied obtaining
\begin{equation}
    \lrr{\Tc_{\bs{r}}+\Kc_{\bs{r}} \Tc^{-1} \lrr{\frac{\Ic}{2} + \Kc}}  \lrr{\eta\bs{J}_L} = \nh\times\bs{E^+}.
    \label{eq:radStP2}
\end{equation}
 \end{color}

To numerically solve \eqref{eq:radStP} \begin{color}{black}and \eqref{eq:radStP2}\end{color}, the discretization scheme will require particular attention.
%a system which can be solved to reconstruct a set of equivalent $\bs{M_L}$ current given the measurement of the $\bs{E}$ field.
\begin{color}{black}Starting with \eqref{eq:radStP}, \end{color}the magnetic current is expanded as
%\begin{equation}
$\bs{M}_L(\bs{r})\approx\sum_{i=1}^{N_e}m_i\bs{f}_i(\bs{r})$ \label{eq:mrwg}
%\end{equation}
where $\{\bs{f}_i\}$ are Rao-Wilton-Glisson (RWG) basis functions (defined here without edge normalization) and $N_e$ is the number of mesh edges. The electric operator $\Tc$ is then tested with rotated RWG functions \cite{cools_accurate_2011} which yields the matrix $\Trwg = \im k \Trwg_s + \im k^{-1} \Trwg_h$, where $[\Trwg_s]_{ij} = \lra{\nh \times \bs{f}_i, \Tc_{s} \bs{f}_j}_{\Gamma}$, $[\Trwg_h]_{ij} = \lra{\nh \times \bs{f}_i, \Tc_{h} \bs{f}_j}_{\Gamma}$,
%\begin{equation}
%    \lrs{\Trwg}_{ij} = \lra{\nh \times \bs{f}_i, \Tc \bs{f}_j}_{\Gamma} \label{eq:Tmat}
%\end{equation}
and $\lra{\bs{a},\bs{b}}_\Gamma=\int_\Gamma\bs a(\bs r) \cdot\bs b(\bs{r}) \diff \bs{r}$.
%\hl{and the same discretization applies to the sub-operators $\Tc_{s}$, $\Tc_h$ for which we define $\Trwg_s$, $\Trwg_h$ respectively}.
%\todo[inline]{That's not correct, Th is discretized by passing derivates, and it also gets a minus sign EC: could this be implementation detail?}
%This obliges the testing basis functions of $\lrr{\Ic/2 + \Kc}$, in the left-following $\lrr{\Ic/2 + \Kc}^{-1}$ term, to be rotated-RWGs which in turn means that the source functions of $\lrr{\Ic/2 + \Kc}$ must be dual elements (we will use in the following Buffa-Christiansen (BC) basis functions, defined in \cite{buffa2007dual, andriulli2008multiplicative}). This follows from the fact that  the identity operator $\Ic$ cannot be effectively  discretized and inverted using only RWGs \cite{andriulli2008multiplicative}.
As a consequence, the magnetic operator must be tested with rotated-RWGs, and to allow for a non-singular discretization of the identity, the source functions used for its discretization must be dual elements \cite{andriulli_multiplicative_2008}---we will use in the following the Buffa-Christiansen (BC) basis functions, a definition of which can be found in \cite{buffa_dual_2007, andriulli_multiplicative_2008}.
%This follows from the fact that  the identity operator $\Ic$ cannot be effectively  discretized and inverted using only RWGs \cite{andriulli2008multiplicative}.
We define the Gram matrix
%\begin{equation}
    $\lrs{\Gbc}_{ij} = \lra{\nh \times \bs{f}_i, \bs{g}_j}_{\Gamma}\,, \label{eq:Gmat}$
%\end{equation}
where $\{\bs g_j(\bs{r})\}$ denote the BC functions
and propose as matrix discretization for the $\Kc$ operator
%\begin{equation}
    $\lrs{\Kbc}_{ij} = \lra{\nh \times \bs{f}_i, \Kc \bs{g}_j}_{\Gamma}\,. \label{eq:Kmat}$
%\end{equation}
Finally, as a consequence if this choice, the source functions of $\Tc_{\bs{r}}$ must be BC functions and a possible choice for the testing functions are rotated-BC basis functions living on $\Gamma_m$.
Thus, we define
%\begin{equation}
%$\lrs{\Tbc_m}_{ij} = \lra{\nh \times \bs{g}_i, \Tc_{\bs{r}} \bs{g}_j}_{\Gamma_m}$,
%\end{equation}
$\Tbc_m = \im k \Tbc_{s,m} + \im k^{-1} \Tbc_{h,m}$ where $\lrs{\Tbc_{s,m}}_{ij} = \lra{\nh \times \bs{g}_i, \Tc_{s,\bs{r}} \bs{g}_j}_{\Gamma_m}$ and $\lrs{\Tbc_{h,m}}_{ij} = \lra{\nh \times \bs{g}_i, \Tc_{h,\bs{r}} \bs{g}_j}_{\Gamma_m}$. From the above choices the discretization of the leftmost $\Kc_{\bs{r}}$ is entirely determined as
%\begin{equation}
    $\lrs{\Krwg_m}_{ij} = \lra{\nh \times \bs{g}_i, \Kc_{\bs{r}} \bs{f}_j}_{\Gamma_m}. \label{eq:Krmat}$
%\end{equation}
%By combining \eqref{eq:Tmat}, \eqref{eq:Gmat}, \eqref{eq:Kmat}, \eqref{eq:Trmat}, and \eqref{eq:Krmat} we get the discretization of the complete operator and the associated equation as
By combining the previous discretization schemes
%\eqref{eq:Tmat}--\eqref{eq:Krmat}
we obtain the discretized equation
\begin{equation}
    \lrr{-\Krwg_m - \Tbc_m\lrr{\Gbc / 2 + \Kbc}^{-1} \Trwg} \lrr{-\Md} = \mbf{e}_m.
    \label{eq:discrete}
\end{equation}
where $\lrs{\mbf e_m}_i=\lra{\nh\times\bs{g}_i,\nh\times\bs{E}^+}_{\Gamma_m}$ and $\mbf{m}$ is the vector of coefficients $m_i$. 
\begin{color}{black}For \eqref{eq:radStP2}, a similar reasoning leads to
\begin{equation}
    \lrr{\Trwg_m + \Kbc_m\Tbc^{-1}\lrr{-\Gbc^\mathrm{T}/2 +\Krwg}} \lrr{\eta\Jd} = \mbf{e}_m\,,
    \label{eq:discrete2}
\end{equation}
with $\Trwg_m = \im k\Trwg_{s,m} + \im k^{-1}\Trwg_{h,m}$, $\lrs{\Trwg_{s,m}}_{ij}=\lra{\nh \times \bs{f}_i, \Tc_{s,\bs{r}} \bs{f}_j}_{\Gamma_m}$, $\lrs{\Trwg_{h,m}}_{ij}=\lra{\nh \times \bs{f}_i, \Tc_{h,\bs{r}} \bs{f}_j}_{\Gamma_m}$, $\lrs{\Kbc_m}_{ij} = \lra{\nh \times \bs{f}_i, \Kc_{\bs{r}} \bs{g}_j}_{\Gamma_m}$, $\Tbc=\im k\Tbc_s + \im k^{-1}\Tbc_h$, $\lrs{\Tbc_s}_{ij} = \lra{\nh \times \bs{g}_i, \Tc_s \bs{g}_j}_{\Gamma}$, $\lrs{\Tbc_h}_{ij} = \lra{\nh \times \bs{g}_i, \Tc_h \bs{g}_j}_{\Gamma}$, $\lrs{\Krwg}_{ij} = \lra{\nh \times \bs{g}_i, \Kc_{\bs{r}} \bs{f}_j}_{\Gamma}$ and $\Jd$ is the vector of coefficients $j_i$ of the electric current expansion $\bs{J}_L(\bs{r})\approx\sum_{i=1}^{N_e}j_i\bs{f}_i(\bs{r})$.

The reader should note that it is not necessary to solve both \eqref{eq:discrete} and  \eqref{eq:discrete2} to obtain both currents: once one of the two currents has been obtained (discretized with RWGs), the discretization of the other as a linear combination of BCs can be obtained after back substitution in \eqref{eq:steklov}. Moreover only one current is required to compute the field in the outside region by using \eqref{eq:radStP} or \eqref{eq:radStP2}, respectively, following the discretization strategies delineated above with the sole difference that the last operators must be evaluated in the point of interest, and not tested with primal or dual functions. 

\section{Quasi-Helmholtz stabilization}
\label{sec:helmoltz}
\begin{figure*}[!t]
\begin{color}{black}
    \centering
    \begin{align}
        \Mbc^{-1}\lrr{\lrr{\Gbc/2 + \Kbc}^{-1} \Trwg} \Mrwg &= \lrr{\sqrt{k}\Psh + \frac{1}{\im\sqrt{k}}\Pl}\lrr{  \lrr{\Gbc/2 + \Kbc}^{-1} \lrr{\im k \Trwg_s + \frac{\im}{k}\Trwg_h}}\lrr{\frac{1}{\sqrt{k}} \Plh + \im\sqrt{k} \Ps} \nonumber\\
         &= \Psh \lrr{\Gbc/2 + \Kbc}^{-1} \lrr{\im k \Trwg_s}\Plh
         + i k \Psh \lrr{\Gbc/2 + \Kbc}^{-1} \lrr{\im k \Trwg_s + \im k^{-1}\Trwg_h}\Ps+\nonumber\\
         &+(i k)^{-1}\Pl \lrr{\Gbc/2 + \Kbc}^{-1} \lrr{\im k \Trwg_s }\Plh
         +\Pl \lrr{\Gbc/2 + \Kbc}^{-1} \lrr{\im k \Trwg_s }\Ps+\mathcal{O}\lrr{k}\label{eq:Stek1}\\  
         &= - \Psh \lrr{\Gbc/2 + \Kbc}^{-1} \Trwg_h\Ps+\Pl \lrr{\Gbc/2 + \Kbc}^{-1} \Trwg_s \Plh+\mathcal{O}\lrr{k}\nonumber\\
        \lrr{\Mbc^{-1}\lrr{\Tbc^{-1}\lrr{-\Gbc^T/2 + \Krwg}} \Mrwg}^{-1}&= \lrr{\sqrt{k}\Plh + \frac{1}{\im\sqrt{k}}\Ps}
        \lrr{  \lrr{-\Gbc^T/2 + \Krwg}^{-1} \lrr{\im k \Tbc_s + \frac{\im}{k}\Tbc_h}}\lrr{\frac{1}{\sqrt{k}} \Psh + \im\sqrt{k} \Pl} \nonumber\\
         &= \Plh \lrr{-\Gbc^T/2 + \Krwg}^{-1} \lrr{\im k \Tbc_s}\Psh
         + i k \Plh \lrr{-\Gbc^T/2 + \Krwg}^{-1} \lrr{\im k \Tbc_s + \im k^{-1}\Tbc_h}\Pl\nonumber\\
         &+(i k)^{-1}\Ps \lrr{\Krwg-\Gbc^T/2}^{-1} \lrr{\im k \Tbc_s }\Psh
         +\Ps \lrr{\Krwg-\Gbc^T/2}^{-1} \lrr{\im k \Tbc_s }\Pl+ \mathcal{O}\lrr{k}\label{eq:Stek2}\\  
         &= - \Plh \lrr{-\Gbc^T/2 + \Krwg}^{-1} \Tbc_h\Pl+\Ps \lrr{-\Gbc^T/2 + \Krwg}^{-1} \Tbc_s \Psh+\mathcal{O}\lrr{k}\nonumber         
    \end{align}
    \hrulefill
    \end{color}
\end{figure*}

%Because of the presence of the EFIE operator, the condition number of the system in \eqref{eq:discrete} increases when decreasing the frequency. This is the well know low-frequency breakdown of the EFIE operator\cite{andriulli2009analysis} \todo[inline]{PR: citation ok?} that here impacts the entire inverse source system.
The linear system  in \eqref{eq:discrete2} inherits the well-known low-frequency breakdown of the EFIE operator, that causes, among other things, the conditioning of the system to grow unbounded as the frequency decreases \cite{eibert_iterative-solver_2004}, \cite{andriulli_well-conditioned_2013}, at the same time the linear system in \eqref{eq:discrete} will behave, frequency-wise, like an MFIE operator requiring low-frequency stabilization \cite{merlini2020magnetic}. Note that some of the standard inverse source formulations in literature may also suffer from similar low-frequency problems and may benefit from what we will propose below. In this contribution however, for the sake of brevity, we will limit the analysis to the low-frequency stabilization of our new formulations only. Define
%\begin{equation}
    % \setcounter{equation}{22}
    $\Mrwg=\Plh k^{-1/2}+\im\Ps k^{1/2}$, $\ \Mbc=\Psh k^{-1/2}+\im\Pl k^{1/2}$,
    \label{eq:precond}
%\end{equation}
where 
%\begin{gather}
    $\Ps=\Sig(\Sig^T\Sig)^+\Sig^T$, $\Plh=\mbf{I}-\Ps$, $\Pl=\Lam(\Lam^T\Lam)^+\Lam^T$, $\Psh=\mbf{I}-\Pl\ $
%\end{gather}
are the quasi-Helmholtz projectors defined respectively in the RWG space and in the dual BC space, $\mbf{I}$ is the identity matrix, and where $\Sig$, $\Lam$, 
%$\mbf{H}$ 
are the star-to-RWG and loop-to-RWG
%, and global-loop-to-RWG 
transformation matrices, the definitions of which can be found in \cite{andriulli_well-conditioned_2013}. We indicate with $+$ the Moore-Penrose (MP) pseudoinverse operation.
%By multiplying the system matrix in \eqref{eq:discrete} with $\Mbc$ on the left and $\Mrwg$ on the right, we obtain \eqref{eq:stable1} where $-\Md=\Mrwg \mbf{x}$.
We propose the following regularization schemes for \eqref{eq:discrete} and \eqref{eq:discrete2}, respectively
\begin{align}
    \Mbc \lrr{-\Krwg_m - \Tbc_m \lrr{\Gbc/2 + \Kbc}^{-1} \Trwg} \Mrwg \mbf{x} &= \Mbc \mbf{e}_m,\label{eq:stable1}\\
    \Mrwg\lrr{\Trwg_m + \Kbc_m\Tbc^{-1}\lrr{-\Gbc^T/2 +\Krwg}}\Mrwg \mbf{y} &= \Mrwg\mbf{e}_m \label{eq:stable2}
   \end{align}
where $-\Md=\Mrwg \mbf{x}$ and $\eta\Jd=\Mrwg \mbf{y}$. The frequency stability of the above equations will now be shown in two steps. First we will show that quasi-Helmholtz projectors can successfully regularize the Steklov-\poin\ operators in both discretizations  presented here. This is proven in \eqref{eq:Stek1} and \eqref{eq:Stek2} where, in addition to standard cancellation properties of projectors on solenoidal spaces \cite{adrian_electromagnetic_2021}, we used: in \eqref{eq:Stek2} the result 
$\| \Ps \lrr{-\Gbc^T/2 + \Krwg}^{-1} \Pl \| =\mathcal{O}\lrr{k^2}$ which follows from  $\| \Ps \lrr{-\Gbc^T/2 + \Krwg} \Pl \| = \mathcal{O}\lrr{k^2}$ (proven in Section~IV.B.1 of \cite{adrian_electromagnetic_2021}) after following a similar procedure as the one in Appendix B of  \cite{adrian_electromagnetic_2021}; in  \eqref{eq:Stek1} the result
 $\|\Pl\lrr{\Gbc/2+\Kbc}^{-1}\Ps\|= \mathcal{O}\lrr{k^2}$ which 
can be proven in a similar and dual way. Finally, the frequency regularity of \eqref{eq:stable1} follows by noticing that $\Mbc\Krwg_m\Mrwg$ is frequency stable \cite{merlini2020magnetic} and that
$\Mbc\Tbc_m \lrr{\Gbc/2 + \Kbc}^{-1} \Trwg\Mrwg=\lrr{\Mbc\Tbc_m\Mbc}\lrr{\Mbc^{-1}  \lrr{\Gbc/2 + \Kbc}^{-1} \Trwg\Mrwg}$ which, following the above developments and the regularity of $\Mbc\Tbc_m\Mbc$, is the product of two frequency regular operators and thus is frequency regular. Dually the  stability and well-conditioning of \eqref{eq:stable2} is proved with $\Mrwg\Kbc_m\Tbc^{-1}\lrr{-\Gbc^T/2 +\Krwg} \Mrwg=\lrr{\Mrwg\Kbc_m\Mbc}\lrr{\Mbc^{-1}\Tbc^{-1}\lrr{-\Gbc^T/2 +\Krwg} \Mrwg}$ and the frequency regularity of $\Mrwg\Kbc_m\Mbc$ (on simply-connected geometries),  $\Mbc^{-1}\Tbc^{-1}\lrr{-\Gbc^T/2 +\Krwg} \Mrwg$, and $\Mrwg\Trwg_m \Mrwg$.

\section{Numerical Results and Discussion}
\label{sec:results}
In the first set of tests we will focus on validating the Steklov-\poin\ approach mapping electric fields to electric fields \eqref{eq:discrete}, a most relevant setting for real case scenarios. 
To test the performance of this formulation, the electric field of an ideal Hertzian dipole at frequency $f=\SI{5}{\giga\hertz}$ is sampled on a spherical surface $1\lambda$ away from a spherical equivalent surface $\Gamma$ of radius $0.67\lambda$ with $\lambda=2\pi/k$.
Our work is compared to the ordinary double current (MP pseudoinversion of the first row of \eqref{eq:linsys})
and non-Love single current (MP pseudoinversion of the upper left block of \eqref{eq:linsys}) formulations. 
The right hand side is the tangential electric field for all formulations.
%In the single current, as for the Steklov-\poin\, the magnetic equivalent current is sought as solution. In all three methods the right hand side is the electric field and the Moore-Penrose pseudoinverse is applied to the discretized matrix. 
The reconstruction errors in $\Omega^+$ are reported in Fig.~\ref{fig:recon} and show comparable performance for all single and double source formulations.
%: there are no major differences among the different formulations which are indeed equivalent in terms of reconstruction capability as one would expect from degrees-of-freedom theory. The reconstruction error increases in the near-field due to the common ill-conditioning of the matrix caused by the evanescent fields.
\begin{figure}
    \centering
    \includegraphics[height=5cm]{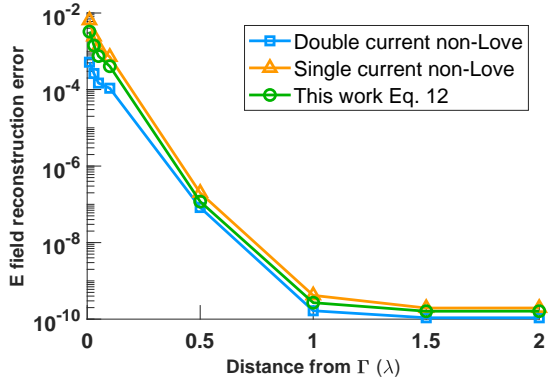}
    \caption{Relative error (spectral norm) of the reconstruction of the electric field for a frequency $f=\SI{5}{\giga\hertz}$ with a spherical reconstruction surface $\Gamma$ of radius $\SI{0.04}{\metre}$ and measurements obtained on a spherical surface of same center as $\Gamma$ and situated \num{1}$\lambda$ away from $\Gamma$.}
    \label{fig:recon}
\end{figure}
To verify the Love condition for our formulation, the tangential component of the electric field is compared with the obtained magnetic currents. The comparison is made on the mesh of $\Gamma$ scanned with $15^\circ\leq\theta\leq165^\circ$ and fixed azimuth angle $\varphi$. Results are shown in Fig. \ref{fig:tan_E} and confirm that for non-Love formulations $\mbf{M}\neq-\nh\times\mbf{E}^+$.
Instead our formulations correctly provides Love currents for which $\mbf{M}_L=-\nh\times\mbf{E}^+$.
As a consequence the internal electric fields radiated by the Love solutions of our formulation is expected to be zero inside the equivalent surface (within the discretization error). This is verified in Fig.~\ref{fig:section} where the magnitude of the electric field is displayed on the plane $z=0$ for different formulations. Clearly  our formulations provides fields that are  orders of magnitude lower than the others inside the equivalent surface $\Gamma$.
\begin{figure}
    \centering
    \includegraphics[height=5cm]{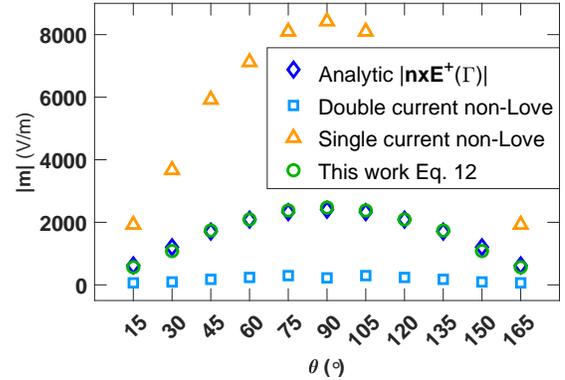}
    \caption{Modulus of the magnetic current on $\Gamma(r)$.}
    \label{fig:tan_E}
\end{figure}
\begin{figure}
    \centering
    \includegraphics[height=6cm]{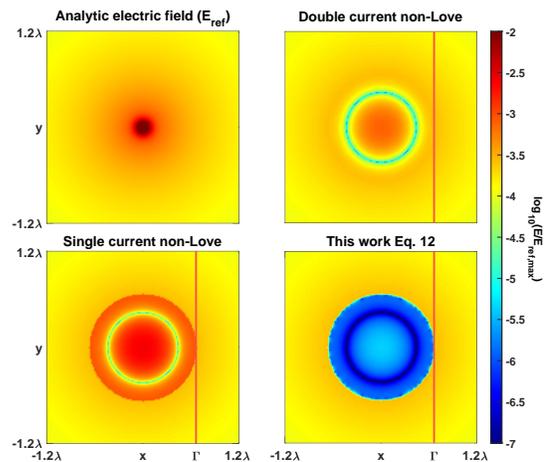}
    \caption{Modulus of the electric field in a $xy$ planar section of $\mbb{R}^3$: fields are normalized on the maximum value of the reference field.}
    \label{fig:section}
\end{figure}
To evaluate the low-frequency behavior of \eqref{eq:stable1} and \eqref{eq:stable2}, we have computed their condition number for different frequencies and compared it with other non-preconditioned formulations (Fig.~\ref{fig:conditioning}).
%we have compared the condition number of three methods while varying the frequency (Fig.~\ref{fig:conditioning}).
In particular, we have computed the conditioning of each matrix by keeping only the first \num{800} singular values of each matrix, this because the overall system is ill-posed and requires pseudoinversion.
%so that the effect of the low frequency breakdown is not covered by the ill-conditioning of the radiation matrices.
As expected, our equations are the only single current Love formulations that show a stable behavior with frequency.
%the double current method suffers from low-frequency breakdown while the single current method has a conditioning, solely dependent on the pseudoinversion cutting point, which is stable. In addition, the Steklov-\poin{} based method has a behavior similar to that of the standard single current method once preconditioned.

\begin{figure}
    \centering
    \includegraphics[height=5cm]{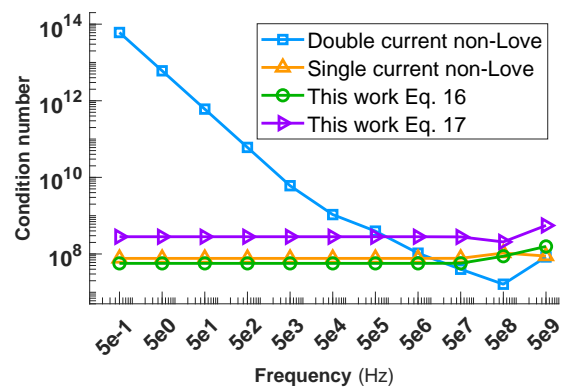}
    \caption{Matrix conditioning in function of the frequency.}
    \label{fig:conditioning}
\end{figure}

\section{Conclusion}
\label{sec:conclusion}

We have presented a new single current approach that naturally yields Love solutions of the inverse source problem
%, without explicit enforcement of additional constraints in the linear system
%that does only require the measurements of one type of field---electric or magnetic.\todo[inline]{EC: I think the information from "that does only require..." is not useful, also all other methods, double and single current methods use one field only generally}
%The method presented in this paper allows to compute the solution of the inverse problem while enforcing the Love condition. Its main and most evident advantage with respect to its alternatives lies in the fact that only the measurements of a single type of field are needed and that it is a single source approach.
%resulting in smaller system matrices. % \todo[inline]{PR/EC: should we not speak about timing at all?}% and consequently allowing the solution of bigger problems.
% However, because of the nature of the operators involved, a correct implementation of this approach needs the introduction of the dual space and of the BC basis functions. In addition, in case both types of current are needed, a further step should be performed after the solution the system to derive the electric currents.
and we have shown that the Love condition is satisfied up to numerical discretization errors and that the tangential components of the fields on the surface are correctly reconstructed. %, which can provide  information for antenna diagnostics comparable to other Love schemes.
The technique is enriched by the first frequency stabilizations of the Steklov-\poin{} operator via quasi-Helmholtz projectors then used to stabilize the new formulation till arbitrary low frequency. This was then confirmed both by theoretical treatments and by numerical results.
   \end{color}

\section*{Acknowledgment}
This work was supported in part by the European Research Council (ERC) under the European Union’s Horizon 2020 research and innovation programme (grant agreement No 724846, project 321), by the European Union’s Horizon 2020 research and innovation programme under the Marie Skłodowska-Curie grant agreement No 955476 (ITN-EID project COMPETE), by the Italian Ministry of University and Research within the Program PRIN2017, EMVISIONING, Grantno. 2017HZJXSZ, CUP:E64I190025300, and by the Italian Ministry of University and Research within the Program FARE, CELER, Grantno. R187PMFXA4.
\bibliographystyle{IEEEtran}
\bibliography{bibliography.bib}

\begin{thebibliography}{10}
\providecommand{\url}[1]{#1}
\csname url@samestyle\endcsname
\providecommand{\newblock}{\relax}
\providecommand{\bibinfo}[2]{#2}
\providecommand{\BIBentrySTDinterwordspacing}{\spaceskip=0pt\relax}
\providecommand{\BIBentryALTinterwordstretchfactor}{4}
\providecommand{\BIBentryALTinterwordspacing}{\spaceskip=\fontdimen2\font plus
\BIBentryALTinterwordstretchfactor\fontdimen3\font minus
  \fontdimen4\font\relax}
\providecommand{\BIBforeignlanguage}[2]{{%
\expandafter\ifx\csname l@#1\endcsname\relax
\typeout{** WARNING: IEEEtran.bst: No hyphenation pattern has been}%
\typeout{** loaded for the language `#1'. Using the pattern for}%
\typeout{** the default language instead.}%
\else
\language=\csname l@#1\endcsname
\fi
#2}}
\providecommand{\BIBdecl}{\relax}
\BIBdecl

\bibitem{lopez_improved_2009}
Y.~A. Lop{\'e}z, F.~L.-H. Andr{\'e}s, M.~R. Pino, and T.~K. Sarkar, ``An
  improved super-resolution source reconstruction method,'' \emph{IEEE
  Transactions on Instrumentation and Measurement}, vol.~58, no.~11, pp.
  3855--3866, 2009.

\bibitem{cappellin_advanced_2011}
E.~J{\o}rgensen, P.~Meincke, and C.~Cappellin, ``Advanced processing of
  measured fields using field reconstruction techniques,'' in \emph{Proceedings
  of the 5th European Conference on Antennas and Propagation (EUCAP)}.\hskip
  1em plus 0.5em minus 0.4em\relax IEEE, 2011, pp. 3880--3884.

\bibitem{foged_practical_2012}
L.~Foged, L.~Scialacqua, F.~Saccardi, J.~A. Quijano, G.~Vecchi, and
  M.~Sabbadini, ``Practical application of the equivalent source method as an
  antenna diagnostics tool [amta corner],'' \emph{IEEE Antennas and Propagation
  Magazine}, vol.~54, no.~5, pp. 243--249, 2012.

\bibitem{petre_planar_1992}
P.~Petre and T.~Sarkar, ``Planar near-field to far-field transformation using
  an equivalent magnetic current approach,'' \emph{IEEE Transactions on
  Antennas and Propagation}, vol.~40, no.~11, pp. 1348--1356, 1992.

\bibitem{sarkar_near-field_1999}
T.~K. Sarkar and A.~Taaghol, ``Near-field to near/far-field transformation for
  arbitrary near-field geometry utilizing an equivalent electric current and
  mom,'' \emph{IEEE Transactions on Antennas and Propagation}, vol.~47, no.~3,
  pp. 566--573, 1999.

\bibitem{quijano_field_2010}
J.~L.~A. Quijano and G.~Vecchi, ``Field and source equivalence in source
  reconstruction on 3d surfaces,'' \emph{Progress In Electromagnetics
  Research}, vol. 103, pp. 67--100, 2010.

\bibitem{kornprobst_solution_2019}
J.~Kornprobst, R.~A.~M. Mauermayer, O.~Neitz, J.~Knapp, and T.~F. Eibert, ``On
  the solution of inverse equivalent surface-source problems,'' \emph{Progress
  In Electromagnetics Research}, vol. 165, pp. 47--65, 2019.

\bibitem{alvarez_reconstruction_2007}
Y.~{\'A}lvarez, F.~Las-Heras, and M.~R. Pino, ``Reconstruction of equivalent
  currents distribution over arbitrary three-dimensional surfaces based on
  integral equation algorithms,'' \emph{IEEE Transactions on Antennas and
  Propagation}, vol.~55, no.~12, pp. 3460--3468, 2007.

\bibitem{eibert_multilevel_2009}
T.~F. Eibert and C.~H. Schmidt, ``Multilevel fast multipole accelerated inverse
  equivalent current method employing rao–wilton–glisson discretization of
  electric and magnetic surface currents,'' \emph{IEEE Transactions on Antennas
  and Propagation}, vol.~57, no.~4, pp. 1178--1185, 2009.

\bibitem{eibert_inverse_2010}
T.~F. Eibert, E.~K. Ismatullah, C.~H. Schmidt \emph{et~al.}, ``Inverse
  equivalent surface current method with hierarchical higher order basis
  functions, full probe correction and multilevel fast multipole
  acceleration,'' \emph{Progress In Electromagnetics Research}, vol. 106, pp.
  377--394, 2010.

\bibitem{jorgensen_improved_2010}
E.~J{\o}rgensen, P.~Meincke, C.~Cappellin, and M.~Sabbadini, ``Improved source
  reconstruction technique for antenna diagnostics,'' in \emph{Proc. 32nd ESA
  Antenna Workshop}, 2010, pp. 1--7.

\bibitem{kornprobst_accuracy_2021}
J.~Kornprobst, J.~Knapp, R.~A.~M. Mauermayer, O.~Neitz, A.~Paulus, and T.~F.
  Eibert, ``Accuracy and {Conditioning} of {Surface}-{Source} {Based}
  {Near}-{Field} to {Far}-{Field} {Transformations},'' \emph{IEEE Transactions
  on Antennas and Propagation}, vol.~69, no.~8, pp. 4894--4908, Aug. 2021.

\bibitem{phaneuf_formulation_2021}
M.~Phaneuf, ``On the formulation and implementation of the {Love}’s
  constraint for the source reconstruction method,'' \emph{IEEE Transactions on
  Antennas and Propagation (Early Access)}, pp. 1--1, 2021.

\bibitem{kornprobst_inverse_2019}
J.~Kornprobst, R.~A. Mauermayer, E.~K{\i}l{\i}{\c{c}}, and T.~F. Eibert, ``An
  inverse equivalent surface current solver with zero-field enforcement by
  left-hand side calder{\'o}n projection,'' in \emph{2019 13th European
  Conference on Antennas and Propagation (EuCAP)}.\hskip 1em plus 0.5em minus
  0.4em\relax IEEE, 2019, pp. 1--3.

\bibitem{eibert_fast_2016}
T.~F. Eibert, D.~Vojvodi{\'c}, and T.~B. Hansen, ``Fast inverse equivalent
  source solutions with directive sources,'' \emph{IEEE Transactions on
  Antennas and Propagation}, vol.~64, no.~11, pp. 4713--4724, 2016.

\bibitem{de_la_bourdonnaye_formulations_1995}
A.~de~La~Bourdonnaye, ``Some formulations coupling finite element and integral
  equation methods for helmholtz equation and electromagnetism,''
  \emph{Numerische Mathematik}, vol.~69, no.~3, pp. 257--268, 1995.

\bibitem{ricci_frequency-stabilized_2022}
P.~Ricci, A.~Merlini, and F.~P. Andriulli, ``On a frequency-stabilized single
  current inverse source formulation,'' \emph{2022 IEEE International Symposium
  on Antennas and Propagation and USNC-URSI Radio Science Meeting (APS/URSI)},
  2022.

\bibitem{cools_accurate_2011}
K.~Cools, F.~Andriulli, D.~De~Zutter, and E.~Michielssen, ``Accurate and
  conforming mixed discretization of the mfie,'' \emph{IEEE antennas and
  wireless propagation letters}, vol.~10, pp. 528--531, 2011.

\bibitem{andriulli_multiplicative_2008}
F.~P. Andriulli, K.~Cools, H.~Bagci, F.~Olyslager, A.~Buffa, S.~Christiansen,
  and E.~Michielssen, ``A multiplicative calderon preconditioner for the
  electric field integral equation,'' \emph{IEEE Transactions on Antennas and
  Propagation}, vol.~56, no.~8, pp. 2398--2412, 2008.

\bibitem{buffa_dual_2007}
A.~Buffa and S.~Christiansen, ``A dual finite element complex on the
  barycentric refinement,'' \emph{Mathematics of Computation}, vol.~76, no.
  260, pp. 1743--1769, 2007.

\bibitem{eibert_iterative-solver_2004}
T.~Eibert, ``Iterative-solver convergence for loop-star and loop-tree
  decompositions in method-of-moments solutions of the electric-field integral
  equation,'' \emph{IEEE Antennas and Propagation Magazine}, vol.~46, no.~3,
  pp. 80--85, 2004.

\bibitem{andriulli_well-conditioned_2013}
F.~P. Andriulli, K.~Cools, I.~Bogaert, and E.~Michielssen, ``On a
  well-conditioned electric field integral operator for multiply connected
  geometries,'' \emph{IEEE transactions on antennas and propagation}, vol.~61,
  no.~4, pp. 2077--2087, 2013.

\bibitem{merlini2020magnetic}
A.~Merlini, Y.~Beghein, K.~Cools, E.~Michielssen, and F.~P. Andriulli,
  ``Magnetic and combined field integral equations based on the quasi-helmholtz
  projectors,'' \emph{IEEE Transactions on Antennas and Propagation}, vol.~68,
  no.~5, pp. 3834--3846, 2020.

\bibitem{adrian_electromagnetic_2021}
S.~B. Adrian, A.~Dely, D.~Consoli, A.~Merlini, and F.~P. Andriulli,
  ``Electromagnetic integral equations: Insights in conditioning and
  preconditioning,'' \emph{IEEE Open Journal of Antennas and Propagation},
  2021.

\end{thebibliography}

\end{document}